\documentclass[12pt]{article}

\usepackage{amsmath, epsfig, cite}
\usepackage{amssymb}
\usepackage{amsfonts}
\usepackage{latexsym}
\usepackage{float}

\newtheorem{thm}{Theorem}[section]

\newtheorem{conj}[thm]{Conjecture}
\newtheorem{lem}[thm]{Lemma}

\numberwithin{equation}{section}

\newcommand{\qed}{{\hfill$\square$}\medskip}

\newcommand*{\pFq}[5]{{}_{#1}F_{#2}\left[ \begin{matrix} #3\\[5pt] #4\end{matrix};#5\right]}

\begin{document}


\begin{center}
{\Large\bf Supercongruences involving $p$-adic Gamma functions}
\end{center}

\vskip 2mm \centerline{Ji-Cai Liu}
\begin{center}
{\footnotesize Department of Mathematics, Wenzhou University, Wenzhou 325035, PR China\\
{\tt jcliu2016@gmail.com} }
\end{center}


\vskip 0.7cm \noindent{\bf Abstract.}
We establish some supercongruences for the truncated ${}_2F_1$ and ${}_3F_2$ hypergeometric series involving the $p$-adic Gamma functions. Some of these results extend the four Rodriguez-Villegas supercongruences on the truncated ${}_3F_2$ hypergeometric series. The corresponding conjectural supercongruences modulo $p^3$ are also proposed for further research.

\vskip 2mm
\noindent{\it MR Subject Classifications}: Primary 11A07, 11S80; Secondary 33C20, 33B15

\section{Introduction}
Rodriguez-Villegas \cite{RV} observed the relationship between the number of points
over $F_p$ on hypergeometric Calabi-Yau manifolds and the truncated hypergeometric series.
Some interesting supercongruences for hypergeometric Calabi-Yau manifolds of dimension $D\le 3$ were conjectured by Rodriguez-Villegas.

To state these results, we first define the truncated hypergeometric series.
For complex numbers $a_i, b_j$ and $z$, with none of the $b_j$ being negative integers or zero, the truncated hypergeometric series are defined as
\begin{align*}
\pFq{r}{s}{a_1,a_2,\cdots,a_r}{b_1,b_2,\cdots,b_s}{z}_n=\sum_{k=0}^{n}\frac{(a_1)_k (a_2)_k\cdots (a_r)_k}{(b_1)_k (b_2)_k\cdots (b_s)_k}\cdot \frac{z^k}{k!},
\end{align*}
where
$(a)_0=1$ and $(a)_k=a(a+1)\cdots (a+k-1)$ for $k\ge 1$.

Throughout this paper let $p\ge 5$ be a prime. Rodriguez-Villegas \cite{RV} proposed four conjectural supercongruences associated to certain modular K3 surfaces.
These were all of the form
\begin{align}
\pFq{3}{2}{\frac{1}{2},&-a,&a+1}{&1,&1}{1}_{p-1}\equiv c_p\pmod{p^2},\label{a1}
\end{align}
where $a=-\frac{1}{2}, -\frac{1}{3},-\frac{1}{4},-\frac{1}{6}$ and $c_p$ is the $p$-th Fourier coefficient of a weight three modular form on a congruence subgroup of $SL(2, \mathbb{Z})$.
The case when $a=-\frac{1}{2}$ was confirmed by van Hamme \cite{vanHamme}, Ishikawa \cite{Ishikawa} and Ahlgren \cite{Ahlgren}.
The other cases when $a=-\frac{1}{3},-\frac{1}{4},-\frac{1}{6}$ were partially proved by
Mortenson \cite{Mortenson-pams}, and finally proved by Z.-W. Sun \cite{Sunzw-aa}.

Let $\langle a\rangle_p$ denote the least non-negative integer $r$ with $a\equiv r \pmod{p}$.
Z.-H. Sun \cite[Theorem 2.5]{Sunzh-jnt} showed that for any $p$-adic integer $a$ with $\langle a\rangle_p\equiv 1\pmod{2}$,
\begin{align}
\pFq{3}{2}{\frac{1}{2},&-a,&a+1}{&1,&1}{1}_{p-1}\equiv 0 \pmod{p^2},\label{a2}
\end{align}
which partially extends \eqref{a1}.
Guo and Zeng \cite[Theorem 1.3]{GZ} have obtained an interesting $q$-analogue of \eqref{a2}.
Using the same idea, Z.-H. Sun \cite[Corollary 2.2]{Sunzh-jnt} also proved that for $\langle a\rangle_p\equiv 1\pmod{2}$,
\begin{align}
\pFq{2}{1}{-a,&a+1}{&1}{\frac{1}{2}}_{p-1}\equiv 0 \pmod{p^2}.\label{a3}
\end{align}
The cases when $a=-\frac{1}{2}, -\frac{1}{3},-\frac{1}{4},-\frac{1}{6}$ on the left-hand side of \eqref{a3} have been dealt with by Z.-H. Sun \cite{Sunzh-pams-2011,Sunzh-jnt}, Z.-W. Sun \cite{Sunzw-scm,Sunzw-ffa} and Tauraso \cite{Tauraso}.

In this paper, we will prove some supercongruences for the truncated ${}_2F_1$ and ${}_3F_2$ hypergeometric series involving the $p$-adic Gamma functions. Some of these results extend the four Rodriguez-Villegas supercongruences on the truncated ${}_3F_2$ hypergeometric series.
Our proof is based on some combinatorial identities involving harmonic numbers and some properties of the $p$-adic Gamma functions.

\begin{thm}\label{t1}
Let $p\ge 5$ be a prime.
For any $p$-adic integer $a$ with $\langle a\rangle_p\equiv 0\pmod{2}$, we have
\begin{align}
\pFq{2}{1}{-a,&a+1}{&1}{\frac{1}{2}}_{p-1}\equiv
(-1)^{\frac{p+1}{2}}\Gamma_p\left(\frac{1}{2}\right)\Gamma_p\left(-\frac{a}{2}\right)\Gamma_p\left(\frac{a+1}{2}\right) \pmod{p^2},\label{a4}
\end{align}
where $\Gamma_p\left(\cdot\right)$ denotes the $p$-adic Gamma function recalled in the next section.
\end{thm}

\begin{thm}\label{t2}
Let $p\ge 5$ be a prime.
For any $p$-adic integer $a$ with $\langle a\rangle_p\equiv 0\pmod{2}$, we have
\begin{align}
\pFq{3}{2}{\frac{1}{2},&-a,&a+1}{&1,&1}{1}_{p-1}\equiv
(-1)^{\frac{p+1}{2}}\Gamma_p\left(-\frac{a}{2}\right)^2\Gamma_p\left(\frac{a+1}{2}\right)^2 \pmod{p^2}.\label{a5}
\end{align}
\end{thm}

In order to prove Theorem \ref{t2}, we need the following supercongruence, which is
a special case of a result due to Z.-H. Sun \cite[Theorem 2.2]{sunzh-a-2015}.

\begin{thm}\label{t3}
Suppose $p\ge 5$ is a prime. For any $p$-adic integer $a$, we have
\begin{align}
\pFq{3}{2}{\frac{1}{2},&-a,&a+1}{&1,&1}{1}_{p-1}\equiv
\pFq{2}{1}{-a,&a+1}{&1}{\frac{1}{2}}_{p-1}^2 \pmod{p^2}.\label{a6}
\end{align}
\end{thm}

Supercongruence \eqref{a6} is a $p$-adic analogue of the following identity:
\begin{align}
\pFq{3}{2}{\frac{1}{2},&-a,&a+1}{&1,&1}{1}=
\pFq{2}{1}{-a,&a+1}{&1}{\frac{1}{2}}^2,\label{a7}
\end{align}
which can be deduced from Clausen's formula. We shall give an alternative proof of
\eqref{a6} by using some combinatorial identities.

The rest of this paper is organized as follows. In the next section, we recall some properties of the $p$-adic Gamma functions, and establish some combinatorial identities involving harmonic numbers.
We prove Theorem \ref{t1} in Section 3, and show Theorems \ref{t2} and \ref{t3} in Section 4.
The corresponding conjectural supercongruences modulo $p^3$ are proposed in the last section.

\section{Some lemmas}
Let $p$ be an odd prime and $\mathbb{Z}_p$ denote the set of all $p$-adic integers. For
$x\in \mathbb{Z}_p$, the Morita's $p$-adic Gamma function \cite[Definition 11.6.5]{Cohen} is defined as
\begin{align*}
\Gamma_p(x)=\lim_{m\to x}(-1)^m\prod_{\substack{0< k < m\\
(k,p)=1}}k,
\end{align*}
where the limit is for $m$ tending to $x$ $p$-adically in $\mathbb{Z}_{\ge 0}$.
We first recall some basic properties of the $p$-adic Gamma function
(see \cite[\S 11.6]{Cohen} for more details on the $p$-adic Gamma function).
For $x\in \mathbb{Z}_p$, we have
\begin{align}
&\Gamma_p(1)=-1,\label{b1}\\
&\Gamma_p(x)\Gamma_p(1-x)=(-1)^{s_p(x)},\label{b2}\\
&\frac{\Gamma_p(x+1)}{\Gamma_p(x)}=
\begin{cases}
-x\quad&\text{if $|x|_p=1$,}\\
-1\quad &\text{if $|x|_p<1$, }
\end{cases}\label{b3}
\end{align}
where $s_p(x)\in \{1,2,\cdots,p\}$ with $s_p(x)\equiv x\pmod{p}$ and
$|\cdot|_p$ denotes the $p$-adic norm.

For $a\in \mathbb{Z}_p$, set $G_1(a)=\Gamma_p^{'}(a)/\Gamma_p(a)$.
We have $G_1(a)\in \mathbb{Z}_p$ (see \cite[Proposition 2.3]{Kilbourn}).
\begin{lem}
Let $p$ be an odd prime. For any $x\in \mathbb{Z}_p$, we have
\begin{align}
G_1(x)\equiv G_1(1)+H_{s_p(x)-1}\pmod{p},\label{b6}
\end{align}
where $H_n$ denotes the $n$-th harmonic number $H_n=\sum_{k=1}^n\frac{1}{k}$.
\end{lem}
{\noindent \it Proof.}
The $p$-adic logarithm is defined as
\begin{align*}
\log_p(1+x)=\sum_{n=1}^{\infty}\frac{(-1)^{n+1}x^n}{n},
\end{align*}
which converges for $x\in \mathbb{C}_p$ with $|x|_p<1$.
Taking the $\log_p$ derivative on both sides of \eqref{b3} gives
\begin{align}
G_1(x+1)-G_1(x)=
\begin{cases}
\frac{1}{x}\quad &\text{if $|x|_p=1$,}\\[5pt]
0\quad &\text{if $|x|_p<1$.}
\end{cases}\label{b7}
\end{align}
For any $p$-adic integer $a$ and $b$ with $a\equiv b \pmod{p}$, by \cite[(2.2) \& (2.3)]{Kilbourn} we have
$\Gamma_p(a)\equiv \Gamma_p(b) \pmod{p}$ and $\Gamma_p^{'}(a)\equiv \Gamma_p^{'}(b) \pmod{p}$,
and so $G_1(a)\equiv G_1(b) \pmod{p}$. Repeatedly applying \eqref{b7}, we obtain
\begin{align*}
G_1(x)&\equiv G_1(s_p(x))\pmod{p}\\
&=G_1(s_p(x)-1)+\frac{1}{s_p(x)-1}\\
&=G_1(1)+H_{{s_p(x)-1}}.
\end{align*}
The desired result is reached.
\qed

We also need some combinatorial identities.
\begin{lem}
For any even integer $n$, we have
\begin{align}
\sum_{k=0}^n{2k\choose k}{n+k\choose 2k}\left(-\frac{1}{2}\right)^k
&=\frac{{n\choose n/2}}{(-4)^{n/2}},\label{b8}\\
\sum_{k=0}^n{2k\choose k}^2{n+k\choose 2k}\left(-\frac{1}{4}\right)^k&=\frac{{n\choose n/2}^2}{4^n},\label{b9}\\
\sum_{k=0}^n{2k\choose k}{n+k\choose 2k}\left(-\frac{1}{2}\right)^k\sum_{i=1}^k\frac{1}{n+i}
&=\frac{{n\choose n/2}}{(-4)^{n/2}}\left(\frac{1}{2}H_n-\frac{1}{2}H_{\frac{n}{2}}\right),\label{b17}\\
\sum_{k=0}^n{2k\choose k}^2{n+k\choose 2k}\left(-\frac{1}{4}\right)^k\sum_{i=1}^k\frac{1}{n+i}
&=\frac{{n\choose n/2}^2}{4^n}\left(\frac{3}{2}H_n-H_{\frac{n}{2}}\right).\label{b18}
\end{align}
\end{lem}

{\noindent \it Proof.}
The identities \eqref{b8} and \eqref{b9} are directly deduced from \eqref{a7} and the following identity (see \cite[(2), page 11]{Bailey}):
\begin{align}
\pFq{2}{1}{-a,&a+1}{&1}{\frac{1}{2}}=\frac{\Gamma\left(\frac{1}{2}\right)}
{\Gamma\left(\frac{1-a}{2}\right)\Gamma\left(1+\frac{a}{2}\right)},\label{b12}
\end{align}
by setting $a=n$.

Note that
\begin{align*}
\sum_{i=1}^k\frac{1}{n+i}=H_{n+k}-H_n.
\end{align*}
In order to prove \eqref{b17} and \eqref{b18}, by \eqref{b8} and \eqref{b9}, it suffices to show that
\begin{align}
\sum_{k=0}^{2n}{2k\choose k}{2n+k\choose 2k}\left(-\frac{1}{2}\right)^k(2H_{2n+k}-3H_{2n}+H_n)&=0,\label{new-1}\\
\sum_{k=0}^{2n}{2k\choose k}^2{2n+k\choose 2k}\left(-\frac{1}{4}\right)^k
(2H_{2n+k}-5H_{2n}+2H_n)&=0.\label{new-2}
\end{align}
Let $A_n$ and $B_n$ denote the left-hand sides of \eqref{new-1} and \eqref{new-2}, respectively.
Using the software package {\tt Sigma} developed by Schneider \cite{schneider-slc-2007}, we find that
$A_n$ and $B_n$ satisfy the following recurrences:
\begin{align*}
(2n+1)A_{n}+2(n+1)A_{n+1}=0,
\end{align*}
and
\begin{align*}
&4 (n+1)^2(2 n+1)^2 (4 n+7) B_n-(4 n+5)(32 n^4+160 n^3+296 n^2+240 n+71)B_{n+1}\notag\\
&+4(n+2)^2 (2 n+3)^2 (4 n+3)B_{n+2}=0,
\end{align*}
respectively.
It is easy to verify that $A_0=0$ and $B_0=B_1=0$, and so $A_n=B_n=0$ for all $n\ge 0$.
\qed

{\noindent \bf Remark.} Combinatorial identities \eqref{b17} and \eqref{b18} can be also automatically discovered and proved by Schneider's computer algebra package {\tt Sigma}.
We refer to \cite[\S 3.1]{schneider-slc-2007} for an interesting approach to finding and proving combinatorial identities of this type.

\section{Proof of Theorem \ref{t1}}
We can rewrite \eqref{a4} as
\begin{align}
\sum_{k=0}^{p-1}{2k\choose k}{a+k\choose 2k}\left(-\frac{1}{2}\right)^k
\equiv
(-1)^{\frac{p+1}{2}}\Gamma_p\left(\frac{1}{2}\right)\Gamma_p\left(-\frac{a}{2}\right)\Gamma_p\left(\frac{a+1}{2}\right) \pmod{p^2}.\label{c1}
\end{align}
Let $\delta$ denote the number $\delta=(a-\langle a\rangle_p)/p$. It is clear that $\delta$ is a $p$-adic integer and $a=\langle a\rangle_p+\delta p$. Since
\begin{align*}
\prod_{i=1}^k(C+x\pm i)=\left(\prod_{i=1}^k(C\pm i)\right)\left(1+x\sum_{i=1}^k\frac{1}{C\pm i}\right)+\mathcal{O}(x^2),
\end{align*}
we have
\begin{align}
{2k\choose k}{a+k\choose 2k}
&={2k\choose k}{\langle a\rangle_p+\delta p+k\choose 2k}\notag\\
&={2k\choose k}\prod_{i=1}^k(\langle a\rangle_p+\delta p+i)\prod_{i=1}^k(\langle a\rangle_p+\delta p+1-i)\prod_{i=1}^{2k}i^{-1}\notag\\
&\equiv
{2k\choose k}{\langle a\rangle_p+k\choose 2k}
\left(1+\delta p\left(\sum_{i=1}^k\frac{1}{\langle a\rangle_p+i}+\sum_{i=1}^k\frac{1 }{\langle a\rangle_p+1-i}\right)\right)\pmod{p^2},\label{c2}
\end{align}
where we have used the fact that ${2k\choose k}\prod_{i=1}^{2k}i^{-1}\in \mathbb{Z}_p$ for $0\le k\le p-1$.
It follows that
\begin{align}
\text{LHS \eqref{c1}}
&\equiv \sum_{k=0}^{p-1}{2k\choose k}{\langle a\rangle_p+k\choose 2k}\left(-\frac{1}{2}\right)^k\notag\\
&\times\left(1+\delta p\left(\sum_{i=1}^k\frac{1}{\langle a\rangle_p+i}+\sum_{i=1}^k\frac{1 }{\langle a\rangle_p+1-i}\right)\right)\pmod{p^2}.\label{c3}
\end{align}

Let $b=p-\langle a\rangle_p$. It is clear that $\langle a\rangle_p\equiv -b\pmod{p}$ and $0\le b-1\le p-1$ is an even integer.  Note that ${-b+k\choose 2k}={b-1+k\choose 2k}$.
Thus,
\begin{align}
&\sum_{k=0}^{p-1}{2k\choose k}{{\langle a\rangle_p+k\choose 2k}}\left(-\frac{1}{2}\right)^k\sum_{i=1}^k\frac{1 }{\langle a\rangle_p+1-i}\notag\\
&\equiv -\sum_{k=0}^{p-1}{2k\choose k}{{-b+k\choose 2k}}\left(-\frac{1}{2}\right)^k\sum_{i=1}^k\frac{1 }{b-1+i}\pmod{p}\notag\\
&=-\sum_{k=0}^{p-1}{2k\choose k}{{b-1+k\choose 2k}}\left(-\frac{1}{2}\right)^k\sum_{i=1}^k\frac{1 }{b-1+i}\notag\\
&\overset{\eqref{b17}}{=}\frac{{b-1\choose (b-1)/2}}{(-4)^{(b-1)/2}}\left(\frac{1}{2}H_{\frac{b-1}{2}}-\frac{1}{2}H_{b-1}\right).\label{c4}
\end{align}
Since ${2n\choose n}/(-4)^n={-1/2\choose n}$ and $b+\langle a\rangle_p=p$, we have
\begin{align}
\frac{{b-1\choose (b-1)/2}}{(-4)^{(b-1)/2}}
={-\frac{1}{2}\choose \frac{b-1}{2}}
\equiv {\frac{p-1}{2}\choose \frac{b-1}{2}}
={\frac{p-1}{2}\choose \frac{\langle a\rangle_p}{2}}
\equiv {-\frac{1}{2}\choose \frac{\langle a\rangle_p}{2}}
=\frac{{\langle a\rangle_p\choose \langle a\rangle_p/2}}{(-4)^{\langle a\rangle_p/2}}\pmod{p}.\label{c5}
\end{align}
It follows from \eqref{c4} and \eqref{c5} that
\begin{align}
&\sum_{k=0}^{p-1}{2k\choose k}{{\langle a\rangle_p+k\choose 2k}}\left(-\frac{1}{2}\right)^k\sum_{i=1}^k\frac{1 }{\langle a\rangle_p+1-i}\notag\\
&\equiv \frac{{\langle a\rangle_p\choose \langle a\rangle_p/2}}{(-4)^{\langle a\rangle_p/2}}
\left(\frac{1}{2}H_{\frac{p-\langle a\rangle_p-1}{2}}-\frac{1}{2}H_{p-\langle a\rangle_p-1}\right)\pmod{p}. \label{c6}
\end{align}

Furthermore, we have
\begin{align}
\sum_{k=0}^{p-1}{2k\choose k}{\langle a\rangle_p+k\choose 2k}\left(-\frac{1}{2}\right)^k
\overset{\eqref{b8}}{=}\frac{{\langle a\rangle_p\choose \langle a\rangle_p/2}}{(-4)^{\langle a\rangle_p/2}},\label{c7}
\end{align}
and
\begin{align}
\sum_{k=0}^{p-1}{2k\choose k}{\langle a\rangle_p+k\choose 2k}\left(-\frac{1}{2}\right)^k\sum_{i=1}^k\frac{1}{\langle a\rangle_p+i}
\overset{\eqref{b17}}{=}\frac{{\langle a\rangle_p\choose \langle a\rangle_p/2}}{(-4)^{\langle a\rangle_p/2}}
\left(\frac{1}{2}H_{\langle a\rangle_p}-\frac{1}{2}H_{\frac{\langle a\rangle_p}{2}}\right).\label{c8}
\end{align}

Finally, combining \eqref{c3} and \eqref{c6}--\eqref{c8} gives
\begin{align}
\text{LHS \eqref{c1}}
\equiv\left(-\frac{1}{4}\right)^{\langle a \rangle_p/2}{\langle a \rangle_p\choose \langle a \rangle_p/2}\left(1+\frac{\delta p}{2}\left(H_{\frac{p-\langle a \rangle_p-1}{2}}-H_{\frac{\langle a \rangle_p}{2}}\right)\right)\pmod{p^2},\label{c9}
\end{align}
where we have used the fact $ H_{p-1-k}\equiv H_k\pmod{p}$ for $0\le k \le p-1$.

Note that
\begin{align}
\frac{\left(\frac{1}{2}\right)_k}{(1)_k}=\frac{{2k\choose k}}{4^k},\label{b14}
\end{align}
and for $a\in \mathbb{Z}_p, n\in \mathbb{N}$
such that none of $a,a+1,\cdots,a+n-1$ in $p\mathbb{Z}_p$
(see \cite[Lemma 17, (4)]{LR}),
\begin{align}
(a)_n=(-1)^n\frac{\Gamma_p(a+n)}{\Gamma_p(a)}.\label{b4}
\end{align}
From \eqref{b14} and \eqref{b4}, we deduce that
\begin{align}
\left(-\frac{1}{4}\right)^{\langle a \rangle_p/2}{\langle a \rangle_p\choose \langle a \rangle_p/2}
\overset{\eqref{b14}}{=}(-1)^{\langle a \rangle_p/2}\frac{\left(\frac{1}{2}\right)_{\langle a \rangle_p/2}}{(1)_{\langle a \rangle_p/2}}
\overset{\eqref{b4}}{=}(-1)^{\langle a \rangle_p/2}\frac{\Gamma_p(1)\Gamma_p\left(\frac{1+\langle a\rangle_p}{2}\right)}{\Gamma_p\left(\frac{1}{2}\right)\Gamma_p\left(1+\frac{\langle a\rangle_p}{2}\right)}.\label{c10}
\end{align}
By \eqref{b2}, we have
\begin{align}
&\Gamma_p\left(\frac{1}{2}\right)^2=(-1)^{\frac{p+1}{2}},\label{c11}\\
&\Gamma_p\left(1+\frac{\langle a\rangle_p}{2}\right)\Gamma_p\left(-\frac{\langle a\rangle_p}{2}\right)=(-1)^{1+\langle a\rangle_p/2}.\label{c12}
\end{align}
Applying \eqref{b1}, \eqref{c11} and \eqref{c12} to the right-hand side of \eqref{c10} and then using
$\langle a\rangle_p=a-\delta p$, we get

\begin{align}
\left(-\frac{1}{4}\right)^{\langle a \rangle_p/2}{\langle a \rangle_p\choose \langle a \rangle_p/2}&=(-1)^{\frac{p+1}{2}}\Gamma_p\left(\frac{1}{2}\right)\Gamma_p\left(\frac{1+\langle a\rangle_p}{2}\right)\Gamma_p\left(-\frac{\langle a\rangle_p}{2}\right)\notag\\
&=(-1)^{\frac{p+1}{2}}\Gamma_p\left(\frac{1}{2}\right)\Gamma_p\left(\frac{1+a-\delta p}{2}\right)\Gamma_p\left(\frac{-a+\delta p}{2}\right).\label{c13}
\end{align}
Note that for $a,b\in\mathbb{Z}_p$ (see \cite[Theorem 14]{LR}),
\begin{align}
\Gamma_p(a+bp)\equiv \Gamma_p(a)\left(1+G_1(a)bp\right) \pmod{p^2}.\label{b5}
\end{align}
Furthermore, applying \eqref{b5} to the right-hand side of \eqref{c13}, we obtain
\begin{align}
\left(-\frac{1}{4}\right)^{\langle a \rangle_p/2}{\langle a \rangle_p\choose \langle a \rangle_p/2}
&\equiv (-1)^{\frac{p+1}{2}}\Gamma_p\left(\frac{1}{2}\right)\Gamma_p\left(\frac{1+a}{2}\right)\Gamma_p\left(-\frac{a}{2}\right)\notag\\
&\times\left(1
+\frac{\delta p}{2}\left(G_1\left(-\frac{a}{2}\right)-G_1\left(\frac{1+a}{2}\right)\right)\right)\pmod{p^2}.\label{c14}
\end{align}
It follows from \eqref{c9} and \eqref{c14} that
\begin{align*}
\text{LHS \eqref{c1}}&\equiv
(-1)^{\frac{p+1}{2}}\Gamma_p\left(\frac{1}{2}\right)\Gamma_p\left(\frac{1+a}{2}\right)\Gamma_p\left(-\frac{a}{2}\right)\\
&\times \left(1
+\frac{\delta p}{2}\left(H_{\frac{p-\langle a \rangle_p-1}{2}}-H_{\frac{\langle a \rangle_p}{2}}+G_1\left(-\frac{a}{2}\right)-G_1\left(\frac{1+a}{2}\right)\right)\right)\pmod{p^2}.
\end{align*}
In order to prove \eqref{c1}, it suffices to show that
\begin{align}
H_{\frac{p-\langle a \rangle_p-1}{2}}-H_{\frac{\langle a \rangle_p}{2}}+G_1\left(-\frac{a}{2}\right)-G_1\left(\frac{1+a}{2}\right)\equiv 0 \pmod{p}.\label{c15}
\end{align}

By \eqref{b6}, we have
\begin{align}
G_1\left(-\frac{a}{2}\right)-G_1\left(\frac{1+a}{2}\right)&\equiv H_{s_p\left(-\frac{a}{2}\right)-1}-H_{s_p\left(\frac{1+a}{2}\right)-1}\pmod{p}.\label{c16}
\end{align}
Since $\langle a\rangle_p$ is an even integer, we have
\begin{align}
s_p\left(-\frac{a}{2}\right)-1&=p-\frac{\langle a\rangle_p}{2}-1,\label{c17}\\
s_p\left(\frac{1+a}{2}\right)-1&=\frac{p+\langle a\rangle_p+1}{2}-1.\label{c18}
\end{align}
Substituting \eqref{c16} into the left-hand side of \eqref{c15} and then using \eqref{c17} and \eqref{c18} gives
\begin{align*}
\text{LHS \eqref{c15}}\equiv
H_{\frac{p-\langle a \rangle_p-1}{2}}-H_{\frac{\langle a \rangle_p}{2}}
+H_{p-\frac{\langle a\rangle_p}{2}-1}-H_{\frac{p+\langle a \rangle_p-1}{2}}\equiv 0 \pmod{p},
\end{align*}
where we have utilized the fact that $H_{p-k-1}\equiv H_k\pmod{p}$ for $0\le k \le p-1$.

\section{Proofs of Theorems \ref{t2} and \ref{t3} }
The proof of Theorem \ref{t2} directly follows from \eqref{a4}, \eqref{a6} and \eqref{c11}. It remains to prove Theorem \ref{t3}. We distinguish two cases to show
\eqref{a6}.

If $\langle a\rangle_p\equiv 1\pmod{2}$, by \eqref{a2} and \eqref{a3}, then \eqref{a6} clearly holds.

If $\langle a\rangle_p\equiv 0\pmod{2}$, by \eqref{c9} and \eqref{b14}, it suffices to show that
\begin{align}
&\sum_{k=0}^{p-1}{2k\choose k}^2{a+k\choose 2k}\left(-\frac{1}{4}\right)^k\notag\\
&\equiv\left(\frac{1}{4}\right)^{\langle a \rangle_p}{\langle a \rangle_p\choose \langle a \rangle_p/2}^2\left(1+\delta p\left(H_{\frac{p-\langle a \rangle_p-1}{2}}-H_{\frac{\langle a \rangle_p}{2}}\right)\right)\pmod{p^2}.\label{d1}
\end{align}

Applying \eqref{c2} to the left-hand side of \eqref{d1} yields
\begin{align}
\text{LHS \eqref{d1}}
&\equiv \sum_{k=0}^{p-1}{2k\choose k}^2{\langle a\rangle_p+k\choose 2k}\left(-\frac{1}{4}\right)^k\notag\\
&\times\left(1+\delta p\left(\sum_{i=1}^k\frac{1}{\langle a\rangle_p+i}+\sum_{i=1}^k\frac{1 }{\langle a\rangle_p+1-i}\right)\right)\pmod{p^2}.\label{d2}
\end{align}
Using the same idea in the previous section and the identities \eqref{b9} and \eqref{b18}, we can show that
\begin{align}
&\sum_{k=0}^{p-1}{2k\choose k}^2{{\langle a\rangle_p+k\choose 2k}}\left(-\frac{1}{4}\right)^k\sum_{i=1}^k\frac{1 }{\langle a\rangle_p+1-i}\notag\\
&\equiv \frac{{\langle a\rangle_p\choose \langle a\rangle_p/2}^2}{4^{\langle a\rangle_p}}
\left(H_{\frac{p-\langle a\rangle_p-1}{2}}-\frac{3}{2}H_{p-\langle a\rangle_p-1}\right)\pmod{p}, \label{d3}\\[10pt]
&\sum_{k=0}^{p-1}{2k\choose k}^2{\langle a\rangle_p+k\choose 2k}\left(-\frac{1}{4}\right)^k
=\frac{{\langle a\rangle_p\choose \langle a\rangle_p/2}^2}{4^{\langle a\rangle_p}},\label{d4}
\end{align}
and
\begin{align}
\sum_{k=0}^{p-1}{2k\choose k}^2{\langle a\rangle_p+k\choose 2k}\left(-\frac{1}{4}\right)^k\sum_{i=1}^k\frac{1}{\langle a\rangle_p+i}
=\frac{{\langle a\rangle_p\choose \langle a\rangle_p/2}^2}{4^{\langle a\rangle_p}}
\left(\frac{3}{2}H_{\langle a\rangle_p}-H_{\frac{\langle a\rangle_p}{2}}\right).\label{d5}
\end{align}
Combining \eqref{d2}--\eqref{d5} gives
\begin{align*}
\text{LHS \eqref{d1}}\equiv \left(\frac{1}{4}\right)^{\langle a \rangle_p}{\langle a \rangle_p\choose \langle a \rangle_p/2}^2\left(1+\delta p\left(H_{\frac{p-\langle a \rangle_p-1}{2}}-H_{\frac{\langle a \rangle_p}{2}}\right)\right)\pmod{p^2}.
\end{align*}
The desired result is reached.

\section{Some open conjectures}
Long and Ramakrishna \cite[Theorem 3]{LR} have extended the case when $a=-\frac{1}{2}$ in \eqref{a1} to a supercongruence modulo $p^3$. Numerical calculation suggests that the other three cases when $a=-\frac{1}{3},-\frac{1}{4},-\frac{1}{6}$ in \eqref{a1} have similar modulo $p^3$ extensions. These three conjectural supercongruences are stated as follows.
\begin{conj}
Let $p\ge 5$ be a prime. The following supercongruences hold modulo $p^3$:
\begin{align}
&\pFq{3}{2}{\frac{1}{2},&\frac{1}{3},&\frac{2}{3}}{&1,&1}{1}_{p-1}
\equiv
\begin{cases}
(-1)^{\frac{p+1}{2}}\Gamma_p\left(\frac{1}{6}\right)^2\Gamma_p\left(\frac{1}{3}\right)^2 \quad&\text{if $p\equiv 1\pmod{6}$},\\[10pt]
(-1)^{\frac{p-1}{2}}\frac{p^2}{18}\Gamma_p\left(\frac{1}{6}\right)^2\Gamma_p\left(\frac{1}{3}\right)^2  \quad&\text{if $p\equiv 5\pmod{6}$},
\end{cases}\label{s1}\\[20pt]
&\pFq{3}{2}{\frac{1}{2},&\frac{1}{4},&\frac{3}{4}}{&1,&1}{1}_{p-1}
\equiv
\begin{cases}
(-1)^{\frac{p+1}{2}}\Gamma_p\left(\frac{1}{8}\right)^2\Gamma_p\left(\frac{3}{8}\right)^2 \quad&\text{if $p\equiv 1,3\pmod{8}$},\\[10pt]
(-1)^{\frac{p-1}{2}}\frac{3p^2}{64}\Gamma_p\left(\frac{1}{8}\right)^2\Gamma_p\left(\frac{3}{8}\right)^2  \quad&\text{if $p\equiv 5,7\pmod{8}$},
\end{cases}\label{s2}\\[20pt]
&\pFq{3}{2}{\frac{1}{2},&\frac{1}{6},&\frac{5}{6}}{&1,&1}{1}_{p-1}
\equiv
\begin{cases}
-\Gamma_p\left(\frac{1}{12}\right)^2\Gamma_p\left(\frac{5}{12}\right)^2 \quad&\text{if $p\equiv 1\pmod{4}$},\\[10pt]
-\frac{5p^2}{144}\Gamma_p\left(\frac{1}{12}\right)^2\Gamma_p\left(\frac{5}{12}\right)^2  \quad&\text{if $p\equiv 3\pmod{4}$}.\label{s3}
\end{cases}
\end{align}
\end{conj}

There is strong numerical evidence to suggest that supercongruence \eqref{a5} also holds modulo $p^3$.
\begin{conj}
Let $p\ge 5$ be a prime.
For any $p$-adic integer $a$ with $\langle a\rangle_p\equiv 0\pmod{2}$, we have
\begin{align}
\pFq{3}{2}{\frac{1}{2},&-a,&a+1}{&1,&1}{1}_{p-1}\equiv
(-1)^{\frac{p+1}{2}}\Gamma_p\left(-\frac{a}{2}\right)^2\Gamma_p\left(\frac{a+1}{2}\right)^2 \pmod{p^3}.\label{s4}
\end{align}
\end{conj}
It is clear that \eqref{s4} reduces to the first cases of \eqref{s1}--\eqref{s3} when $a=-\frac{1}{3},-\frac{1}{4},-\frac{1}{6}$.

Unfortunately, the method in this paper is not applicable for proving these conjectures.

\vskip 5mm \noindent{\bf Acknowledgments.} 
The author is grateful to Professor Zhi-Hong Sun for helpful comments on this manuscript.
The author also would like to thank the anonymous referee for careful reading of this manuscript and valuable suggestions, which make the paper more readable.


\begin{thebibliography}{99}
\small \setlength{\itemsep}{-.8mm}

\bibitem{Ahlgren}S. Ahlgren, Gaussian hypergeometric series and combinatorial congruences, Symbolic computation, number theory, special functions, physics and combinatorics (Gainesville, FL, 1999), 1--12, Dev. Math., 4, Kluwer, Dordrecht, 2001.

\bibitem{Bailey}W.N. Bailey, Generalized Hypergeometric Series, Stechert-Hafner, 1964.

\bibitem{Cohen}H. Cohen, Number Theory. Vol. II. Analytic and Modern Tools, Grad. Texts in Math., vol. 240, Springer, New York, 2007.

\bibitem{GZ}V.J.W. Guo and J. Zeng, Some $q$-supercongruences for truncated basic hypergeometric series, Acta Arith. 171 (2015), 309--326.

\bibitem{Ishikawa}T. Ishikawa, Super congruence for the Ap\'ery numbers, Nagoya Math. J. 118 (1990), 195--202.

\bibitem{Kilbourn}T. Kilbourn, An extension of the Ap\'ery number supercongruence, Acta Arith. 123 (2006), 335--348.

\bibitem{LR}L. Long and R. Ramakrishna, Some supercongruences occurring in truncated hypergeometric series, Adv. Math. 290 (2016), 773--808.

\bibitem{Mortenson-pams}E. Mortenson, Supercongruences for truncated ${}_{n+1}F_n$ hypergeometric series with applications to certain weight three newforms, Proc. Amer. Math. Soc. 133 (2005), 321--330.

\bibitem{RV}F. Rodriguez-Villegas, Hypergeometric families of Calabi-Yau manifolds, Calabi-Yau varieties and mirror symmetry (Toronto, ON, 2001), Fields Inst. Commun., vol. 38, Amer. Math. Soc., Providence, RI, 2003, 223--231.

\bibitem{schneider-slc-2007}C. Schneider, Symbolic summation assists combinatorics, S\'em. Lothar. Combin. 56 (2007), B56b, 36 pp.

\bibitem{Sunzh-pams-2011}Z.-H. Sun, Congruences concerning Legendre polynomials, Proc. Amer. Math. Soc. 139 (2011), 1915--1929.

\bibitem{Sunzh-jnt}Z.-H. Sun, Generalized Legendre polynomials and related supercongruences, J. Number Theory 143 (2014), 293--319.
    
\bibitem{sunzh-a-2015}Z.-H. Sun, Note on super congruences modulo $p^2$, preprint, 2015, 	arXiv:1503.03418.

\bibitem{Sunzw-scm}Z.-W. Sun, Super congruences and Euler numbers, Sci. China Math. 54 (2011), 2509--2535.

\bibitem{Sunzw-aa}Z.-W. Sun, On sums involving products of three binomial coefficients, Acta Arith. 156 (2012), 123--141.

\bibitem{Sunzw-ffa}Z.-W. Sun, Supercongruences involving products of two binomial coefficients, Finite Fields Appl. 22 (2013), 24--44.

\bibitem{Tauraso}R. Tauraso, An elementary proof of a Rodriguez-Villegas supercongruence, preprint, 2009, arXiv:0911.4261.

\bibitem{vanHamme}L. van Hamme, Proof of a conjecture of Beukers on Ap\'ery numbers, Proceedings of the conference on $p$-adic analysis (Houthalen, 1987), 189--195.

\end{thebibliography}
\end{document}